\title{Families of trees decompose the random graph in any arbitrary way}

\author{Raphael Yuster
\thanks{
e-mail: raphy@research.haifa.ac.il \qquad
World Wide Web: http:$\backslash\backslash $research.haifa.ac.il$\backslash$\symbol{126}raphy}
\\ Department of Mathematics\\ University of
Haifa at Oranim\\ Tivon 36006, Israel}

\date{} 

\documentstyle [11pt]{article}
\newtheorem{theo}{Theorem}[section]

\newtheorem{lemma}[theo]{Lemma}
\newtheorem{coro}[theo]{Corollary}

\newcommand\npf{\mbox{ }\hfill\sqr\vskip6pt}
\def\sqr{$\vcenter{\hrule height.2mm
\hbox{\vrule width.2mm height2mm\kern2mm
\vrule width.2mm}\hrule height.2mm}$}
\newcommand{\ignore}[1]{}

\begin{document}
\maketitle
\setcounter{page}{1}
\begin{abstract}
Let $F=\{H_1,\ldots,H_k\}$ be a family of graphs. A graph $G$ with $m$ edges is called
{\em totally $F$-decomposable} if for {\em every} linear combination of the form
$\alpha_1 e(H_1) + \cdots + \alpha_k e(H_k) = m$ where each $\alpha_i$ is a nonnegative integer,
there is a coloring of the edges of $G$ with $\alpha_1+\cdots+\alpha_k$ colors such that
exactly $\alpha_i$ color classes induce each a copy of $H_i$, for $i=1,\ldots,k$.
We prove that if $F$ is any fixed family of trees then 
$\log n/n$ is a sharp threshold function for the property that the random graph $G(n,p)$ is
totally $F$-decomposable.
In particular, if $H$ is a tree, then $\log n/n$ is a sharp threshold function for the property
that $G(n,p)$ contains $\lfloor e(G)/e(H) \rfloor$ edge-disjoint copies of $H$.
\end{abstract}

\section{Introduction}
All graphs considered here are finite, undirected and have no loops or multiple edges.
For the standard terminology used the reader is referred to \cite{Bo}.
For the standard terminology used in Random Graph Theory  the reader is referred to \cite{Bo2}.
Let $H$ and $G$ be two graphs. An $H$-packing of $G$ is a collection of edge-disjoint subgraphs of
$G$, each being isomorphic to $H$. The {\em $H$-packing number of $G$}, denoted $P(H,G)$,
is the maximum size of an $H$-packing of $G$. Clearly, $P(H,G) \leq \lfloor e(G)/e(H) \rfloor$.
If equality holds, we say that $G$ has an {\em optimal $H$-packing}. If, in addition,
$e(H)$ divides $e(G)$ and $P(H,G)=e(G)/e(H)$ then we say that $G$ has an
{\em $H$-decomposition}.

Packing and decomposition theory is a central topic in Graph Theory and Design Theory.
We shall mention here the following general results. If $G=K_n$, and $n$ is sufficiently large,
Wilson \cite{Wi} gave necessary and sufficient conditions for the existence of an
$H$-decomposition of $K_n$. For graphs $H$ with at most 5 vertices, necessary and sufficient
conditions for an $H$-decomposition are known for all $n$ (cf. \cite{CoDi}).
Caro and Yuster \cite{CaYu} gave a closed formula for $P(H,K_n)$, for $n$ sufficiently large.
The formula only depends on the degree sequence of $H$, and on $n$. It follows that for $n$
sufficiently large, there are necessary and sufficient conditions for the existence of an optimal
$H$-packing. For arbitrary graphs $G$, and for fixed graphs $H$ other than trees, almost
nothing is known. We mention the result of Gustavsson \cite{Gu}, who gave necessary and
sufficient conditions for an $H$-decomposition of a graph $G$, where $G$ is ``almost'' complete.
In case $H$ is a tree, it has been proved \cite{Yu2} that every graph $G$ with minimum degree
$d(G) \geq \lfloor n/2 \rfloor$ ($n$ sufficiently large), has an optimal $H$-packing. This result
is sharp for all trees with at least two edges. A related result concerning trees, appearing in \cite{Yu1},
shows that every graph $G$ which is a good expander has an optimal $H$-packing. By ``good expander''
we mean that $G$ is $\Omega(\sqrt{n \log n})$ edge-expanding, and $n$ is sufficiently large as
a function of $H$.

Let $F=\{H_1,\ldots,H_k\}$ be a family of graphs. A graph $G$ with $m$ edges is called
{\em totally $F$-decomposable} if for {\em every} linear combination of the form
$\alpha_1 e(H_1) + \cdots + \alpha_k e(H_k) = m$ where each $\alpha_i$ is a nonnegative integer,
there is a coloring of the edges of $G$ with $\alpha_1+\cdots+\alpha_k$ colors such that
exactly $\alpha_i$ color classes induce each a copy of $H_i$, for $i=1,\ldots,k$.
In other words, $G$ is totally $F$-decomposable if we can decompose it into elements of $F$
in any arbitrary way. Clearly, if $H$ is a graph then, considering the very special case of the family
$F=\{H,K_2\}$, we have that $G$ is totally $F$-decomposable if and only if $G$ has an optimal
$H$-packing.

Let $G=G(n,p)$ denote, as usual, the random graph with $n$ vertices 
and edge probability $p$. In the extensive study of the properties
of random graphs, many researchers observed that
there are sharp {\em threshold functions} for various natural 
graph properties.
For a graph property $A$ and for a function $p=p(n)$, we say that
$G(n,p)$ satisfies $A$ {\em almost surely} if the probability that
$G(n,p(n))$ satisfies $A$ tends to $1$ as $n$ tends to infinity.
We say that a function $f(n)$ is a {\em sharp threshold function}
for the property $A$ if there are two positive constants $c$ and
$C$ so that $G(n,cf(n))$ almost surely does not satisfy $A$ and
$G(n,Cf(n))$ satisfies $A$ almost surely. The seminal result of Friedgut and Kalai \cite{FrKa}
states that every monotone (w.r.t. edge addition) graph property has a sharp threshold
function.

In this paper we consider the property of being totally $F$-decomposable
where $F$ is any fixed family of trees.
The property of being totally $F$-decomposable is not monotone. In fact, even the very special
case of the property of having an optimal $H$-packing is not monotone
for every tree $H$ with at least three edges.
Let $m$ be a positive integer, and let $G$ be any graph having
$e(H) \cdot m - 1$ edges, and having an optimal $H$-packing.  Add two isolated
vertices to $G$, and denote the new graph by $G'$. $G'$ also has an
optimal $H$-packing. Now add to $G'$ an edge between the two isolated vertices.
The new graph has $e(H) \cdot m$ edges, but, obviously, does not have
an $H$-decomposition. One may claim that the non-connectivity of $G'$
caused the non-monotonicity. However, it is not difficult to show that if $H$ has
three edges or more, there exist connected graphs $G$ with an optimal
$H$-packing, and such that it is possible to add an edge to $G$ and obtain
a graph which does not have an optimal $H$-packing. For example, let $H=K_{1,k}$
where $k \geq 3$ and let $G=K_{1,sk-1}$ where $s \geq 1$. $G$ contains $s-1$
edge-disjoint copies of $H$ and therefore has an optimal $H$-packing. Add to $G$
an edge connecting two nonroot vertices. The new graph has $sk$ edges but does
not contain $s$ edge-disjoint copies of $H$.

It is easy to show
that if $p=\sqrt{\frac{\log n}{n}}$ then $G(n,p)$ almost surely is
$\Theta(\sqrt{n \log n})$ edge-expanding, and thus, by the result in \cite{Yu1},
$G(n,p)$ almost surely has an optimal $H$-packing. On the other hand,
it is well-known that if $p=c\log n/n$ where $c$ is a sufficiently small constant,
then $G(n,p)$ has isolated vertices and many small components. Thus,
trivially, we almost surely do not have $\lfloor e(G)/e(H) \rfloor$ edge-disjoint
copies of $H$ in $G$, for any fixed tree $H$ with at least two edges.
Consequently, {\em if there exists} a sharp threshold function $p(n)$ for the property of
containing an optimal $H$-packing then it must be within these bounds.
In this paper we prove that, indeed, such a sharp threshold function exists.
In fact, we show something much stronger:
\begin{theo}
\label{t11}
Let $F$ be a family of trees. Then, $p(n)= \frac{\log n}{n}$ is
a sharp threshold function for the property of being totally $F$-decomposable.
\end{theo}
By considering $F=\{H,K_2\}$ we have the following immediate corollaries:
\begin{coro}
\label{c12}
Let $H$ be a fixed tree with at least two edges. Then, $p(n)= \frac{\log n}{n}$ is
a sharp threshold function for the property of having an optimal $H$-packing.
\end{coro}
\begin{theo}
\label{c13}
Let $H$ be a fixed tree with at least two edges. Then, there are absolute positive
constants $c$ and $C$ such that if $p \leq c \frac{\log}{n}$ then  $G(n,p)$ almost surely does not have
an $H$-decomposition, and if $p \geq C\frac{\log n}{n}$ then $G(n,p)$ has an $H$-decomposition with
probability approaching $1/e(H)$, as $n \rightarrow \infty$. 
\end{theo}
 
The next section contains the proof of Theorem \ref{t11}.
The final section contains some concluding remarks.

\section{Proof of the main result}
Let $F=\{H_1,\ldots,H_k\}$ be a family of trees, and let $h_i=e(H_i)$ denote the number of edges of
$H_i$. Put $c=h_1+\ldots+h_k$. Clearly, we can assume  $c \geq 2$, otherwise there is nothing to prove.
Let $C\geq (28c)^{30}$. We prove the following
\begin{lemma}
\label{l21}
Let $p=C\frac{\log n}{n}$. Then, $G(n,p)$ is almost always surely totally $F$-decomposable.
\end{lemma}
The proof of Theorem \ref{t11} follows immediately from Lemma \ref{l21} and the trivial fact
that for a sufficiently small constant $c$, $G(n,c\log n/n)$ is almost surely not totally $F$-decomposable.
We note here that the constant $(28c)^{30}$ can easily be improved. We make no attempt to
optimize it.

In the rest of this paper we assume $n$ is sufficiently large, whenever necessary.
The first part of our proof does not concern random graphs. We show that if $G$ is an $n$-vertex
graph, that has several ``semi-random'' properties (to be stated in the following lemma), then $G$ is totally
$F$-decomposable.
Since $G(n,p)$ will almost surely have these semi-random properties, the result will follow.

\begin{lemma}
\label{l22}
Let $G$ be a graph with $n$ vertices and assume that the edges of $G$ are colored red and blue.
Let $G_1$ be the spanning subgraph consisting of the red edges and let
$G_2$ be the spanning subgraph consisting of the blue edges.
Furthermore, assume that
\begin{enumerate}
\item
$\Delta(G) \leq 1.5C \log n$.
\item
$\delta(G_1) \geq 0.05C\log n$.
\item
$\delta(G_2) \geq 0.5C\log n$.
\item
Every subset $X$ of vertices with $|X| \leq n/2$
is incident with at least $0.46C|X|\log n$ blue edges whose other endpoint is not in $X$.
\item
For every subset $X$ of vertices with $|X| \geq n/2$
there are at least $0.05Cn \log n$ blue edges with both endpoints in $X$.
\end{enumerate}
Then, $G$ is totally $F$-decomposable.
\end{lemma}

The proof of Lemma \ref{l22} is based, in part, on the following lemma,
which establishes similar ``semi-random'' conditions for being $H$-decomposable.
\begin{lemma}
\label{l23}
Let $H$ be a tree with $h \geq 2$ edges. Let $C_1\geq(10h)^{10}$.
Suppose $G$ is an $n$-vertex graph with $mh$ edges where $m$ is an integer.
Furthermore, suppose that
\begin{enumerate}
\item
$\Delta(G) \leq 1.5C_1 \log n \qquad \delta(G) \geq 0.4C_1 \log n$.
\item
Every subset $X$ of vertices with $|X| \leq n/2$
is incident with at least $0.42C_1|X|\log n$ edges whose other endpoint is not in $X$.
\end{enumerate}
Then, $G$ has an $H$-decomposition.
\end{lemma}

The proof of Lemma \ref{l23} is very similar to the proof of the main result appearing
in \cite{Yu1}. We say that a graph $G=(V,E)$ is
$r$ {\em edge-expanding} if for every nonempty $X \subset V$ with $|X| \leq
|V|/2$, there are at least $r|X|$ edges between $X$ and $V \setminus X$.
The following result is proved in \cite{Yu1}.
\begin{lemma}
\label{l24}
Let $H$ be any tree with $h \geq 1$ edges. Let $G$ be a graph on $n$ vertices
and $mh$ edges where $m$ is an integer.
If $G$ is $10h^4\sqrt{n \log n}$ edge-expanding
then $G$ has an $H$-decomposition.
\end{lemma}
Notice the difference between Lemma \ref{l24} and Lemma \ref{l23}.
Lemma \ref{l24} requires a large edge expansion, namely $\Theta(\sqrt{n \log n})$,
but there are no constraints placed on the maximum degree.
Theorem \ref{l23} proves a seemingly stronger result, where the expansion needed is
only $\Theta(\log n)$ but there is also a maximum-degree and minimum-degree constraint.
It is not difficult to modify the (rather complicated) proof of Lemma \ref{l24} to use a smaller edge expansion
under the assumption that the maximum and minimum degrees are bounded as in Lemma
\ref{l23}. For completeness, this modified proof yielding Lemma \ref{l23} appears in Appendix A.

\noindent
{\bf Proof of Lemma \ref{l22}:}\,
Fix an $n$-vertex graph $G=(V,E)$ and a two edge-coloring satisfying the conditions of Lemma \ref{l22}.
Let $\alpha_1,\ldots,\alpha_k$ be nonnegative integers satisfying
$\alpha_1h_1+\ldots+\alpha_k h_k=|E|$. We must show how to decompose
$G$ into $\alpha_i$ copies of $H_i$ for $i=1,\ldots,k$.

We will partition $F$ into two parts $F_1,F_2$ as follows. If $\alpha_i < \frac{|E|}{20c^2}$
then $H_i \in F_1$, otherwise $H_i \in F_2$. Note that it is possible that $F_1 = \emptyset$, but,
obviously, we must always have $F_2 \neq \emptyset$. Put $|F_1|=s$ and $|F_2|=k-s$, and assume,
w.l.o.g. that $F_1 = \{H_1,\ldots,H_s\}$.

\noindent
We first take care of the $\alpha_i$ copies of $H_i$ for $H_i \in F_1$. All these copies will be taken
from $G_1$ (thus, using only red edges), by the following greedy process.
Assume we have already picked some $r < \alpha_1+\ldots+\alpha_s$ copies from $G_1$
and we wish to pick the next tree which should be isomorphic to some specific $H_i \in F_1$.
The $r$ previously
selected trees use less than $\alpha_1 h_1+ \cdots +\alpha_s h_s < \frac{|E|}{20c}$ edges.
Notice that condition 1 in Lemma \ref{l22} implies $|E| \leq 0.75Cn \log n$ and condition 2 implies that
the number of red edges is at least $0.025Cn \log n$. Thus, the number of red edges is at least $|E|/30$.
Since $|E|/30-|E|/(20c) \geq |E|/120$ there is a subgraph on the remaining unused red edges
with minimum degree at least $|E|/(120n) \gg h_i$. Thus, there is a copy of $H_i$ using
previously unused red edges.

\noindent
Let $G'=(V,E')$ denote the spanning subgraph of $G$ composed of the blue edges and the remaining red
edges not used in the copies of $H_1,\ldots,H_s$. Notice that $|E'|=\alpha_{s+1}h_{s+1}+\ldots+\alpha_k h_k$.
It remains to show how $G'$ has a decomposition into $\alpha_i$ copies of $H_i$ for $i=s+1,\ldots,k$.
For $i=s+1,\ldots,k$ define:
$$
t_i = \left\lfloor 2000c^2 \frac{\alpha_i}{\alpha_{s+1}+\ldots+\alpha_k} \right\rfloor.
$$
We have $t_i \geq 100$ since $\alpha_i \geq |E|/(20c^2)$ and since, clearly,
$\alpha_{s+1}+ \ldots +\alpha_k < |E|$. Thus,
\begin{equation}
\label{e11}
100 \leq t_i \leq 2000c^2.
\end{equation}
Given any set of trees, we can {\em concatenate} them
into one tree by choosing one vertex from each tree, and identifying all the chosen vertices. The
concatenated tree is, by definition, decomposable to its originators. Let $H$ denote the tree
obtained by concatenating $t_i$ copies of $H_i$ for each $i=s+1,\ldots,k$. Note that
$H$ has exactly $h=t_{s+1}h_{s+1}+\ldots+t_k h_k$ edges. By (\ref{e11}), $h \leq 2000c^3$.
Now define:
$$
q = \left\lfloor 0.99 \frac{|E'|}{h} \right\rfloor.
$$

\noindent
{\bf Claim:}\, $t_i q \leq \alpha_i$ for $i=s+1, \ldots, k$.\\
{\bf Proof:}\, It suffices to prove that:
$$
0.99 \frac{|E'|}{h} \cdot 2000c^2 \frac{\alpha_i}{\alpha_{s+1}+\ldots+\alpha_k} \leq \alpha_i.
$$
or equivalently
\begin{equation}
\label{e12}
0.99 \frac{\alpha_{s+1}h_{s+1} + \ldots + \alpha_k h_k}{h} 2000c^2 \leq \alpha_{s+1}+\ldots+\alpha_k.
\end{equation}
We will use the fact that
$$
t_i = \left\lfloor 2000c^2 \frac{\alpha_i}{\alpha_{s+1}+\ldots+\alpha_k} \right\rfloor \geq
2000c^2 \frac{\alpha_i}{\alpha_{s+1}+\ldots+\alpha_k} -1 \geq 1980c^2 \frac{\alpha_i}{\alpha_{s+1}+
\ldots+\alpha_k}.
$$
We therefore have
$$
h=t_{s+1}h_{s+1} + \ldots +t_k h_k \geq  1980c^2 \frac{\alpha_{s+1}h_{s+1} + \ldots + \alpha_k h_k}
{\alpha_{s+1}+ \ldots + \alpha_k},
$$
and, therefore, (\ref{e12}) holds. This completes the proof of the claim.

\noindent
According to the last claim, we can define $b_i=\alpha_i-t_iq$ for $i=s+1,\ldots,k$ and we are guaranteed
that the $b_i$ are nonnegative integers. Our next goal is to find in $G'$ a spanning subgraph $G''$ with
the property that $G''$ has a decomposition in which there are exactly $b_i$ copies of $H_i$ for each
$i=s+1,\ldots,k$, and $\Delta(G'') \leq 0.04C \log n$. We use the following procedure. Assume that we have
already found a subgraph $\hat{G}$ of $G'$ with $\Delta(\hat{G}) \leq 0.04C \log n$ and which contains a decomposition into $b_i$ copies
of each $H_i$, $i=s+1,\ldots,k-1$ and $b_k-1$ copies of $H_k$ (completing the last element is, clearly, the
most difficult situation in the construction, as we may assume $h_k$ is the size of the largest tree in $F_2$).
We wish to add a copy of $H_k$ to $\hat{G}$ such that the edges of $H_k$ are taken from $E' \setminus 
e(\hat{G})$, and such that the resulting graph $G''$ has $\Delta(G'') \leq 0.04C\log n$. We first estimate the number of
edges in $\hat{G}$:
\begin{equation}
\label{e13}
e(\hat{G}) < h_{s+1}b_{s+1}+ \ldots + h_k b_k = \sum_{i=s+1}^k h_i(\alpha_i-t_iq) = |E'| - qh \leq
\end{equation}
$$
|E'|-  h\left(0.99\frac{|E'|}{h}-1\right) =0.01|E'|+ h \leq 0.01|E'|+2000c^3.
$$
It follows that $\hat{G}$ has at least $\lceil n/2 \rceil$ vertices with degrees not exceeding
$(0.04|E'|+8000c^3)/n$. Let $X$ be such a set of $\lceil n/2 \rceil$ vertices. Consider the graph
induced by the vertices of $X$ and the edges of $E' \setminus e(\hat{G})$. We denote this graph by
$X$ as well. By the fifth condition in Lemma \ref{l22}, we have that
the number of edges of $X$ is at least $0.05Cn \log n-0.01|E'|-2000c^3 >0.01Cn\log n$.
Hence $X$ has a subgraph with minimum degree $\Theta(\log n) \gg h_k$ and
thus, we can find in $X$ a copy of $H_k$.
Joining the edges of a copy of $H_k$ in
$X$ to $\hat{G}$ we obtain the graph $G''$ which, by construction, is a subgraph of $G'$ and, furthermore,
$$
\Delta(G'') \leq \max\left\{\Delta(\hat{G}) \; , \; \frac{0.04|E'|+8000c^3}{n}+\Delta(H_k)\right\} \leq
$$
$$
\max\left\{0.04C\log n \; , \; 0.03C\log n+\frac{8000c^3}{n}+h\right\} = 0.04C\log n.
$$
Having constructed the graph $G''$ we now come to the final stage of the proof. Denote by
$G^*$ the spanning subgraph of $G'$ obtained by deleting the edges of $G''$. We claim that
$G^*$ has an  $H$-decomposition, and the number of elements in this decomposition is $q$.
We prove this using Lemma \ref{l23}. First, we must show that $e(G^*)=q \cdot e(H)=qh$.
This is true since:
$$
e(G^*)=|E'|-e(G'')=\sum_{i=s+1}^k \alpha_i h_i - \sum_{i=s+1}^k b_ih_i =
q \sum_{i=s+1}^k h_it_i = qh.
$$
Next, we show that $G^*$ and $H$ satisfy the other conditions of Lemma \ref{l23} with $m=q$
and $C_1=C$.
First notice that
$$
C_1=C \geq (28c)^{30} \geq (10h)^{10}.
$$
Thus, $C_1$ satisfies the assumption in Lemma \ref{l23}. Also,
$$
\Delta(G^*) \leq \Delta(G) \leq 1.5C\log n = 1.5C_1 \log n.
$$
$$
\delta(G^*) \geq \delta(G') - \Delta(G'') \geq \delta(G_2) - \Delta(G'') \geq 0.5C \log n - 0.04C \log n > 0.4C_1\log n.
$$
So Condition 1 of Lemma \ref{l23} is fulfilled.
Finally by Condition 4 in Lemma \ref{l22}, and by the fact that $\Delta(G'') \leq 0.04C\log n$ we have that
every subset of vertices $X$ with $|X| \leq n/2$ has at least
$$
0.46C|X| \log n - 0.04C|X|\log n \geq 0.42C|X| \log n
$$
edges connecting $X$ and $V \setminus X$ in $G^*$. So Condition 2 of Lemma \ref{l23} is fulfilled,
and it follows from lemma \ref{l23} that $G^*$ has an $H$-decomposition into $q$ copies.

Since every copy of $H$ is decomposable into
$t_i$ copies of $H_i$ for each $i=s+1,\ldots,k$, we have that $G^*$ has a decomposition into
$qt_i$ copies of $H_i$ for each $i=s+1,\ldots,t$.
\noindent
It is now easy to see that $G'$ has a decomposition into $\alpha_i$ copies of $H_i$ for
$i=s+1,\ldots,k$. First note that by our construction, $G''$ and $G^*$ are edge disjoint
and their edges union is $E'$. Finally notice that the decomposition of
$G''$ has $b_i$ copies of $H_i$ and the decomposition of $G^*$ has $qt_i$ copies of $H_i$. Together,
this gives $b_i+qt_i=\alpha_i$ copies of $H_i$.
\npf

Finally, we prove that $G(n,p)$ almost surely satisfies the properties stated in Lemma \ref{l22}.
\begin{lemma}
\label{l25}
Let $p=C\frac{\log n}{n}$.
Almost surely, the edges of $G(n,p)$ can be colored using the two colors red and blue such that
if $G_1$ is the spanning subgraph of the red edges and $G_2$ is the spanning subgraph of the blue
edges then
\begin{enumerate}
\item
$\Delta(G) \leq 1.5C \log n$.
\item
$\delta(G_1) \geq 0.05C\log n$.
\item
$\delta(G_2) \geq 0.5C\log n$.
\item
Every subset $X$ of vertices with $|X| \leq n/2$
is incident with at least $0.46C|X|\log n$ blue edges whose other endpoint is not in $X$.
\item
For every subset $X$ of vertices with $|X| \geq n/2$
there are at least $0.05Cn \log n$ blue edges with both endpoints in $X$.
\end{enumerate}
\end{lemma}
{\bf Proof:}\,
Each pair of vertices $(i,j)$ will be a red edge with probability $p/15$, a blue edge with probability $14p/15$
and a non-edge with probability $1-p$. The resulting graph is a typical element of $G(n,p)$ where the
edges have two possible colors.
The expected degree of a vertex $v$ in $G_1$, is $C\log n/15$.
By standard large deviation estimates (cf. \cite{AlSp} Appendix A),
\begin{eqnarray}
\nonumber \Pr \left[|d_{G_1}(v)- \frac{1}{15}C\log n| > \frac{1}{60}C\log n\right] < \\
\nonumber 2\exp\left(-\frac{(C\log n/60)^2}{2C\log n/15}+\frac{(C\log n/60)^3}{2(C\log n/15)^2}\right) = \\
\nonumber 2\exp\left(-\frac{1}{640}C\log n\right) < \frac{1}{n^2}.
\end{eqnarray}
Thus, with probability at least $1-1/n$ the second part of the lemma holds.
The expected degree of a vertex $v$ in $G_2$, is $14C\log n/15$. As in the previous inequality,
$$
\Pr \left[|d_{G_2}(v)- \frac{14}{15}C\log n| > \frac{6.5}{15}C\log n\right] <
2\exp\left(-\frac{169}{3136}C \log n\right) < \frac{1}{n^2}.
$$
Thus, with probability at least $1-1/n$ the third part of the lemma holds.
Notice that since $d_{G_1}(v)+d_{G_2}(v)=d_G(v)$ the previous two inequalities show that
$d_G(v) \leq 1.45C\log n$ for all $v \in V$, and in particular the first part of the lemma holds with
probability at least $1-2/n$.

\noindent
Let $V$ denote the set of all vertices, and consider a nonempty $X \subset V$ with
$|X| \leq n/2$. Let $out(X)$ denote the number of blue edges connecting
a vertex of $X$ and a vertex of $V\setminus X$ in $G_2$. The expectation of $out(X)$
is $\frac{14}{15}p|X|(n-|X|) \geq \frac{7}{15}C|X| \log n$. Using large deviation once again we get:
$$
\Pr \left[out(X) - \frac{14}{15}p|X|(n-|X|) < - 0.01p|X|(n-|X|)\right] < \exp\left(-\frac{0.01^2}{2\frac{14}{15}}p|X|(n-|X|)\right) <
\frac{1}{n^2{n \choose {|X|}}}.
$$
Since there are ${n \choose {|X|}}$ subsets of size $|X|$, and since there are $n/2$ sizes
to consider, we have that with probability at least $1-1/n$, for every $X \subset V$ with
$|X| \leq n/2$,
$$
out(X) > \frac{14}{15}p|X|(n-|X|)-0.01p|X|(n-|X|) > 0.92p|X|(n-|X|) \geq 0.46C|X|\log n,
$$
which means that the fourth part of the lemma holds.

\noindent
Finally, consider $X \subset V$ with
$|X| \geq n/2$. Let $in(X)$ denote the number of blue edges with both endpoints in $X$.
The expectation of $in(X)$ is $\frac{14}{15}p\frac{|X|(|X|-1)}{2}$. Hence,
$$
\Pr \left[in(X) - \frac{14}{15}p\frac{|X|(|X|-1)}{2} < - \frac{7}{15}p\frac{|X|(|X|-1)}{2} \right] <
\exp\left(-\frac{(\frac{7}{15})^2}{2\frac{14}{15}}p\frac{|X|(|X|-1)}{2} \right) <
\frac{1}{n^2{n \choose {|X|}}}.
$$
Since there are ${n \choose {|X|}}$ subsets of size $|X|$, and since there are $n/2$ sizes
to consider, we have that with probability at least $1-1/n$, for every $X \subset V$ with
$|X| \geq n/2$,
$$
in(X) \geq \frac{14}{15}p\frac{|X|(|X|-1)}{2}-\frac{7}{15}p\frac{|X|(|X|-1)}{2}=
\frac{7}{15}p\frac{|X|(|X|-1)}{2} > \frac{7}{15}p\frac{n^2}{9} > 0.05Cn\log n
$$
which means that the fifth part of the lemma holds. \npf

\noindent
Now, lemma \ref{l21} clearly follows from Lemma \ref{l22} and Lemma \ref{l25}. \npf

\section{Concluding remarks}
\begin{itemize}
\item
The proof of Lemma \ref{l21} is algorithmic. Namely, given as input a graph taken from the
probability distribution $G(n,p)$, and given $\alpha_1,\ldots,a_k$ such that
$\alpha_1h_1+\ldots+\alpha_kh_k=e(G)$, the algorithm {\em almost surely} finds
a decomposition of $G$ into $\alpha_i$ copies of $H_i$ for $i=1,\ldots,k$. This follows from the
fact, proved in Lemma \ref{l25}, that the input graph almost surely satisfies the conditions of
Lemma \ref{l22}, and from the fact that all the details in Lemma \ref{l22},
except for the final part which uses Lemma \ref{l23},
can be implemented by a {\em deterministic} polynomial time algorithm. Indeed, we only need to count degrees,
and to find trees in graphs whose minimum degree is higher than the number of vertices of the tree. These
computational tasks are easy to perform in polynomial time. Finally, Lemma \ref{l23} has a randomized polynomial
time algorithm (exactly as shown for Lemma \ref{l24} in \cite{Yu1}).
\item
It may be interesting to find other families of graphs for which sharp threshold functions can be determined
for the property of being totally decomposable. More specifically, say $Q$ is the four vertex graph consisting
of a triangle and an additional edge. Let $F=\{Q,K_2\}$. Can one determine a sharp threshold function for
being $F$ decomposable? Notice that for trivial divisibility reasons, some families do not have an associated
nontrivial threshold function. For example, suppose  $C_3 \in F$. Every graph $G=(V,E)$ which is totally
$F$-decomposable must either have $|E| \neq 0 \bmod 3$ or else have all its degrees even.
For every nontrivial $p$, $G(n,p)$ does not satisfy these two requirements with probability very close to $1/3$.
\item
A proof similar to that of Lemma \ref{l22}, combined with the main result of \cite{Yu1}
yields the following theorem, whose proof is omitted.
\begin{theo}
\label{t31} 
Let $F$ be a family of trees. Then, for $n$ sufficiently large, every graph with minimum degree
$0.5n(1+o_n(1))$ is totally $F$-decomposable.
\end{theo}
Notice that for the very special case of $F=\{H,K_2\}$, the result in \cite{Yu2}
states that Theorem \ref{t31} is true even without the error term $o_n(1)$.
It would be interesting to determine if this stronger version of Theorem \ref{t31} holds
for every fixed family of trees.
\end{itemize}

\newpage
\appendix
\section{Proof of Lemma \ref{l23}}
For the rest of this appendix we assume $G=(V,E)$ is an $n$-vertex graph
with $mh$ edges, where $m$ is an integer and
which satisfies the conditions of Lemma \ref{l23}. Namely, $\delta(G) \geq 0.4C_1 \log n$,
$\Delta(G) \leq 1.5C_1\log n$, and $out(X) \geq 0.42C_1|X|\log n$ for every $X \subset V$ with
$|X| \leq n/2$, where $out(X)$ is the number of edges between $X$ and $V \setminus X$.
Notice that $m \geq \frac{C_1}{5h}n\log n$ and recall that $C_1 \geq (10h)^{10}$.
\begin{lemma}
\label{la1}
$E$ can be partitioned into $h$ subsets $E_1,\ldots,E_h$, each having size $m$, such that each
of the spanning subgraphs $G_i=(V,E_i)$ of $G$, and each vertex $v$ have the following properties:
\begin{enumerate}
\item
$$
\left|d_i(v) - \frac{d(v)}{h}\right| \leq 0.05\frac{d(v)}{h^2},
$$
where $d_i(v)$ denotes the degree of $v$ in $G_i$.
\item
$out_i(X) > \frac{1}{5h}C_1|X| \log n$, for every $X \subset V$ with $|X| \leq n/2$.
Here $out_i(X)$ denotes the number of edges between $X$ and $V\setminus X$ in $G_i$.
\end{enumerate}
\end{lemma}
{\bf Proof:}\,
Each edge of $E$ chooses a random integer between $0$ and $h$, where $0$ is chosen with
probability $\beta=n^{-1/2}$ and the other numbers are chosen with probability
$\alpha=(1-\beta)/h$. All the choices are independent. For $i=0,\ldots,h$ let
$F_i$ denote the set of edges which selected $i$. Let $d'_i(v)$ denote the number of edges
of $F_i$ incident with $v$. The expectation of $|F_i|$ is $\alpha|E|=m(1-\beta)$, for $i\neq 0$.
Using a large deviation inequality we get that for $i\neq 0$:
\begin{equation}
\label{e1}
\Pr [|F_i| > m] = {\rm Prob}[|F_i| - m(1-\beta) > m\beta] <
\exp\left(-\frac{2m^2\beta^2}{mh}\right) =
\end{equation}
$$
\exp\left(-\frac{2m}{nh}\right)\leq \exp\left(-\frac{2\frac{C_1}{5h} n \log n}{nh}\right) < \frac{1}{n}.
$$
For all $i=1,\ldots,h$ and for all $v \in V$ we have:
\begin{equation}
\label{e2}
\Pr \left[\left|d'_i(v) - \alpha d(v)\right| > 0.02\frac{d(v)}{h^2}\right] <
2\exp\left(-\frac{0.0008\frac{d(v)^2}{h^4}}{d(v)}\right) <
2\exp\left(-\frac{0.0003C_1 \log n}{h^4}\right) < \frac{1}{n^2}.
\end{equation}
Similarly, for $i=0$ we have:
\begin{equation}
\label{e3}
\Pr [|d'_0(v) - \beta d(v)| > 0.02\frac{d(v)}{h^2}] < \frac{1}{n^2}.
\end{equation}
From equation (\ref{e1}), (\ref{e2}) and (\ref{e3}) we get that with probability at least
$1-h/n-hn/n^2-n/n^2 > 0.9$ all of the following events happen simultaneously:
\begin{enumerate}
\item
$|F_i| \leq m$ for $i=1,\ldots,h$.
\item
$|d'_i(v) - \alpha d(v)| \leq 0.02\frac{d(v)}{h^2}$ for all $i=1,\ldots,h$ and for all $v \in V$.
\item
$|d'_0(v) - \beta d(v)| \leq 0.02\frac{d(v)}{h^2}$ for all $v \in V$.
\end{enumerate}
Consider, therefore, a partition of $E$ into $F_0, \ldots, F_h$ in which
all of these events hold. Since $|F_i| \leq m$, we may partition $F_0$ into
$h$ subsets $Q_1,\ldots,Q_h$, where $|Q_i|=m-|F_i|$. Put $E_i = F_i
\cup Q_i$ for $i=1,\ldots,h$. Note that $|E_i|=m$ and $E_i \cap E_j =
\emptyset$ for $1 \leq i < j \leq h$. Put $G_i=(V,E_i)$ and let $d_i(v)$ be
the degree of $v$ in $G_i$. Clearly,
\begin{equation}
\label{e4}
d_i(v) \geq d'_i(v) \geq \alpha d(v) - 0.02\frac{d(v)}{h^2} = \frac{d(v)}{h}
- \frac{d(v)}{\sqrt n h} -  0.02\frac{d(v)}{h^2} \geq \frac{d(v)}{h} - 0.03\frac{d(v)}{h^2}.
\end{equation}
We also need to bound $d_i(v)$ from above:
$$
d_i(v) \leq d'_i(v) + d'_0(v) \leq \alpha d(v) + \beta d(v) + 0.04\frac{d(v)}{h^2} =
$$
\begin{equation}
\label{e5}
\frac{d(v)}{h} - \frac{d(v)}{\sqrt n h} + 0.04\frac{d(v)}{h^2} + \frac{d(v)}{\sqrt n} \leq \frac{d(v)}{h}+
0.05\frac{d(v)}{h^2}.
\end{equation}
It now follows from inequalities (\ref{e4}) and (\ref{e5}) that $|d_i(v) -
\frac{d(v)}h| \leq 0.05\frac{d(v)}{h^2}$.

\noindent
It remains to show that with probability greater than $1-0.9=0.1$, the requirements regarding
$out_i(X)$ are met for each $i=1,\ldots,h$ and each $X \subset V$ with $|X| \leq n/2$.
Since $E_i \supset F_i$ it suffices to show that for each such $X$, the number of edges
between $X$ and $V \setminus X$ in the subgraph induced by $F_i$, denoted $out'_i(X)$,
is at least $\frac{1}{5h}C_1|X| \log n$. Consider a subset $X$. The expectation of $out'_i(X)$
is $\alpha \cdot out(X)$. Applying large deviation we get:
$$
\Pr [|out'_i(X)-\alpha \cdot out(X) | > \alpha \cdot out(X)/2] <
2\exp(-\frac{2 \cdot out(X)^2\alpha^2/4}{out(X)})=2\exp(-out(X)\alpha^2/2) \leq
$$
$$
2\exp(-out(X)/(2h^3)) \leq 2\exp(-0.42C_1|X|\log n/(2h^3)) \ll \frac{1}{nh{n \choose |X|}}.
$$
Since there are ${n \choose |X|}$ sets of size $|X|$, and since there are $n/2$ possible
sizes to consider, we get from the last inequality that with probability at
least $0.5 > 0.1$, for all $i=1,\ldots,h$ and for all sets $X \subset V$
with $|X| \leq n/2$,
$$
|out'_i(X) - \alpha \cdot out(X)| \leq \alpha \cdot out(X)/2.
$$
In particular this means that
$$
out_i(X) \geq out'_i(X) \geq \alpha \cdot out(X)/2 \geq \frac{1-1/\sqrt{n}}{h}0.21C_1|X|\log n
\geq \frac{1}{5h}C_1|X|\log n.
$$
\npf

We call a partition of $E$ into the subsets $E_i$ having the properties
guaranteed by Lemma \ref{la1} a {\em feasible partition}. Given a feasible
partition, our next goal is to orient the edges of every
$E_i$, such that the oriented sets, denoted by $E_i^*$ have certain
properties. Let $d^+_i(v)$ and $d^-_i(v)$ denote the outdegree and indegree
of $v$ in $E_i^*$, respectively. Clearly, $d_i(v)=d^+_i(v)+d^-_i(v)$ for all
$v \in V$ and $i=1,\ldots,h$. In order to define the properties which we
require from our orientation, we need several definitions.

Let $q$ be a leaf of $H$. Fix a rooted orientation $H(q)$ of $H$ where the
root of $H$ is $q$. Such an orientation can be obtained by performing a
sequential search of the tree, like Breadth-First Search or Depth First Search.
Let $e_1, \ldots, e_h$ be the oriented edges of $H(q)$, in the order
they are discovered by the search. Note that for $i=2, \ldots, h$, the edge
$e_i=(x,y)$ has a unique {\em parent-edge}, which is the unique edge $e_j$
entering $x$. (Thus, $e_j=(z,x)$ for some $z$). The edge $e_1$ is the only
edge which has no parent, since it is the only edge emanating from $q$. For
$i=2,\ldots,h$, let $p(i)=j$ if $e_j$ is the parent of $e_i$. Note that
$p(i) < i$. We say that $j$ is a {\em descendent} of $i$ if $j=i$ or if
$p(j)$ is a descendent of $i$. Note that this definition is recursive.

An orientation of a feasible partition is called a {\em feasible orientation}
if for all $v \in V$, $d^-_{p(i)}(v) = d^+_i(v)$, where $i=2,\ldots,h$, and
$|d^+_i(v)-d^-_i(v)| \leq i \cdot 0.1d(v)/h^2$, for all $i=1, \ldots,h$. Note
that the second requirement implies also that $|d^+_i(v) - d_i(v)/2| \leq
0.05d(v)/h$ and, similarly, $|d^-_i(v) - d_i(v)/2| \leq 0.05d(v)/h$.

\begin{lemma}
\label{la2}
Every feasible partition has a feasible orientation.
Furthermore, in every feasible orientation
\begin{equation}
\label{e6}
d^+_i(v) \geq 0.17\frac{C_1}{h}\log n
\end{equation}
holds for all $v \in V$ and for all $i=1,\ldots,h$.
\end{lemma}
{\bf Proof:}\, We show how to construct our orientation in $h$ stages,
where in stage $i$ we orient the edges of $E_i$ and form $E_i^*$. We begin by
orienting $E_1$. It is well-known by Euler's Theorem (cf. \cite{Bo}), that
the edges of every undirected graph can be oriented such that the indegree
and outdegree of every vertex differ by at most 1. Such an orientation is
called Eulerian. We therefore let $E_1^*$ be any Eulerian orientation of
$E_1$. Thus $|d^+_1(v)-d^-_1(v)| \leq 1 \leq 0.1 d(v)/h^2$. Assume now
that we have oriented all the subsets $E_j$ for $1 \leq j < i$, such that the
conditions of a feasible orientation hold for $j$. We show how to orient the
edges of $E_i$, such that the conditions also hold for $i$. Let $j=p(i)$, and
put $c_v=d^-_j(v)$. We are required to orient the edges of $E_i$ such that
for every $v \in V$, $d^+_i(v)=c_v$. Our initial goal is to show that
$|d^+_i(v)-d^-_i(v)| \leq i \cdot 0.1d(v)/h^2$. Our second goal is to show that
such an orientation exists. The following inequality achieves the first goal:
$$
|d^+_i(v)-d^-_i(v)|=|2c_v-d_i(v)|=|2d_j(v)-2d^+_j(v)-d_i(v)| \leq
|2d^+_j(v)-d_j(v)|+|d_j(v)-d_i(v)| =
$$
$$
|d^+_j(v)-d^-_j(v)|+|d_j(v)-d_i(v)| \leq j \cdot 0.1 \frac{d(v)}{h^2} + |d_j(v) -
\frac{d(v)}{h}| + |d_i(v) - \frac{d(v)}{h}| \leq
$$
$$
j \cdot 0.1 \frac{d(v)}{h^2} + 0.1\frac{d(v)}{h^2} \leq i \cdot 0.1 \frac{d(v)}{h^2}.
$$

\noindent
We now need to show that the desired orientation exists. Note that $\sum_{v
\in V}c_v=m$ and hence the desired orientation exists if every vertex $v$ can
select $c_v$ edges from the $d_i(v)$ edges adjacent to $v$, and such that
every edge of $E_i$ is selected by exactly one of its endpoints. To prove
this is possible we define a bipartite graph $B$ as follows. $B$ has two
vertex classes of size $m$ each. One vertex class is $E_i$, while the other
vertex class, denoted by $S$, contains $c_v$ copies of each $v$. Thus,
$S=\{v^{(k)} \; | \; v \in V, \; 1 \leq k \leq c_v\}$. The edges of $B$ are
defined as follows. A member $v^{(k)} \in S$ is connected to $e \in E_i$ if $v$
is an endpoint of $e$. Clearly, our aim is to show that $B$ has a perfect
matching. By Hall's Theorem (cf. \cite{Bo}), it suffices to show that for
every set $S' \subset S$, $|N(S')| \geq |S'|$ where $N(S') \subset E_i$ are
the neighbors of $S'$ in $B$. Fix $\emptyset \neq S' \subset S$. Let $V'=\{v
\in V \; | \; v^{(k)} \in S'\}$. Put $V'=\{v_1,\ldots,v_t\}$. Clearly, $|S'|
\leq \sum_{l=1}^t c_{v_l}$. Note that $N(S')$ contains all the edges of $E_i$
which have an endpoint in $V'$. Let $T_1 \subset E_i$ be the set of edges
having only one endpoint in $V'$ and let $T_2=N(S') \setminus T_1$ be the set
of edges of $E_i$ having both endpoints in $V'$. Put $t_1=|T_1|$ and
$t_2=|T_2|$. Clearly, $t_1+2t_2 = \sum_{l=1}^t d_i(v_l)$. We first consider
the case $t \leq n/2$. By Lemma \ref{la1} we have
$out_i(V')=t_1 \geq \frac{1}{5h}C_1t \log n$. Therefore,
$$
|N(S')| = t_1 + t_2 = \sum_{l=1}^t \frac{d_i(v_l)}{2} + \frac{t_1}{2} \geq
\sum_{l=1}^t \frac{d_i(v_l)}{2} + \frac{1}{10h}C_1t\log n >
$$
$$
\sum_{l=1}^t \left(\frac{d_i(v_l)}{2} + \frac{1.5C_1\log n}{15h}\right) \geq
\sum_{l=1}^t \left(\frac{d_i(v_l)}{2} + \frac{d(v_l)}{15h}\right) \geq
\sum_{l=1}^t c_{v_l} \geq |S'|.
$$
The case where $t > n/2$ is proved as follows. Put $V''=V \setminus V' =
\{v_{t+1}, \ldots, v_n\}$. Note that $T_1$ is the set of edges connecting
$V'$ with $V''$. Since
$|V''| \leq n/2$ we have $t_1 \geq \frac{1}{5h}C_1(n-t) \log n$. Now,
$$
|N(S')| = t_1 + t_2 = \sum_{l=1}^t \frac{d_i(v_l)}{2} + \frac{t_1}{2} \geq
\sum_{l=1}^t \frac{d_i(v_l)}{2} + \frac{1}{10h}C_1(n-t)\log n >
$$
$$
m - \sum_{l=t+1}^n \left(\frac{d_i(v_l)}{2} - \frac{1.5C_1\log n}{15h}\right) \geq m -
\sum_{l=t+1}^n c_{v_l} = \sum_{l=1}^t c_{v_l} \geq |S'|.
$$

\noindent
Finally, we need to show that (\ref{e6}) holds. We use the fact that
$|d^+_i(v) - d_i(v)/2| \leq 0.05d(v)/h$ and Lemma \ref{la1} which
states that $|d_i(v) - d(v)/h| \leq 0.05d(v)/h^2$
and the fact that $h \geq 2$ to obtain that
$$
|d^+_i(v) - \frac{d(v)}{2h}| \leq 0.0625\frac{d(v)}{h}.
$$
Thus,
$$
d^+_i(v) \geq \frac{d(v)}{2h} - 0.0625\frac{d(v)}{h} =0.4375\frac{d(v)}{h} \geq 0.175\frac{C_1}{h}\log n.
$$
\npf

A feasible orientation defines a decomposition of the edges of $G$ into a set $L^*$ of $m$
edge-disjoint connected graphs, each graph having $h$ edges, one from each
$E_i$. Furthermore, each of these graphs is homomorphic to $H(q)$ (and, thus,
to $H$), in the sense that every member of $L^*$ which happens to be a tree,
is {\em isomorphic} to $H$. Unfortunately, not all the members of $L^*$ are
necessarily trees, and we will need to mend $L^*$ in order to obtain our
desired decomposition.

We now describe the process which creates $L^*$. Fix a feasible orientation
of $G$, and let $D^+_i(v)
\subset E_i^*$ denote those edges of $E_i^*$ which emanate from $v$, and let
$D^-_i(v) \subset E_i^*$ be the edges of $E_i^*$ which enter $v$. For $i=2,
\ldots, h$ and for all $v \in V$ we know that
$|D^-_{p(i)}(v)|=|D^+_i(v)|=d^+_i(v)$. Therefore, let $B_i(v)$ be a {\em
perfect matching} between $D^-_{p(i)}(v)$ and $D^+_i(v)$. (Note that there
are $d^+_i(v)!$ different ways to select $B_i(v)$, so we pick one
arbitrarily). The members of $B_i(v)$ are, therefore, pairs of edges in the
form $((x,v),(v,y))$ where $(x,v) \in D^-_{p(i)}(v)$ and $(v,y) \in
D^+_i(v)$. We say that $(x,v)$ and $(v,y)$ are {\em matched} if
$((x,v),(v,y)) \in B_i(v)$ for some $i$. The transitive closure of the
"matched" relation defines an equivalence relation where the equivalence
classes are connected directed graphs, each having $h$ edges, one from each
$E_i^*$, and which are homomorphic to $H(q)$, by the homomorphism which maps
the edge $e_i$ of $H(q)$ to the edge belonging to $E^*_i$ in an equivalence
class. Thus, $L^*$ is the set of all of these graphs, (or, in set theoretical
language, the quotient set of the equivalence relation). Note that although
each $T \in L^*$ is homomorphic to $H(q)$, it is not necessarily isomorphic
to $H(q)$ since $T$ may contain cycles. For a simple example, consider the
case where $H(q)$ is a directed path on 3 edges $(q,a,b,c)$. It may be the
case that $T$ is composed of the edges $(x,y) \in E^*_1$, $(y,z) \in E^*_2$
and $(z,x) \in E^*_3$. Thus $T$ is a directed triangle, but not a directed
path on 3 edges. It is clear, however, that if $T$ happens to be a tree, (or,
equivalently, if $T$ contains $h+1$ vertices) then it {\em is} isomorphic to
$H(q)$.

As noted, there are many ways to create $L^*$. In fact, there are
$$
\Pi_{i=2}^h \Pi_{v \in V} d^+_i(v)!
$$
different ways to create the decomposition $L^*$. Our goal is to show that in
at least one of these decompositions, all the members of $L^*$ are, in fact,
trees. Before proceeding with the proof we require a few definitions.

For a member $T \in L^*$, and for $i=1,\ldots,h$, let $T_i$ be the subgraph
of $T$ which consists only of the first $i$ edges, namely those belonging to
$E^*_1 \cup \ldots \cup E^*_i$. Note that $T_i$ is a connected subgraph of
$T$. Let $T(i)$ be the edge of $T$ which belongs to $E^*_i$. Note that for $i
> 1$, $T_i$ is obtained from $T_{i-1}$ by adding the edge $T(i)$. Now,
suppose $T_{i-1}$ is a tree, and $T_i$ is not a tree. Let $T(i)=(v,u)$. (Note
that $(v,u) \in D^+_i(v)$ in this case). It follows that $u$ already appears
in $T_{i-1}$. We therefore call an edge $T(i)=(v,u)$ {\em bad} if $u$ already
appears in $T_{i-1}$. Otherwise, the edge is called {\em good}. Clearly, $T$
is a tree iff all its $h$ edges are good. For $1 \leq i \leq j \leq h$, let
$$
N(v,i,j) = \{T \in L^* \; | T(i) \in D^+_i(v), T(j) \;\; is \;\; bad\}.
$$
Clearly, $|N(v,i,j)| \leq d^+_i(v)$. Our next goal is to show that if all the
$n(h-1)$ perfect matchings $B_i(v)$ are selected randomly and independently,
then with high probability, $|N(v,i,j)|$ is significantly smaller than
$d^+_i(v)$.

\begin{lemma}
\label{la3}
If all the perfect matchings $B_i(v)$ are selected randomly and
independently, then with probability at least $0.9$, for all $i=1,\ldots,h$,
for all $j=i,\ldots,h$ and for all $v \in V$, $|N(v,i,j)| \leq \log n$.
\end{lemma}
{\bf Proof:}\,
Since the perfect matchings are selected
randomly and independently, we may assume that the $n$ matchings $B_j(u)$ for
all $u \in V$ are selected {\em after} all the other $n(h-2)$ matchings
$B_k(u)$, for $k \neq j$, are selected. Prior to the selection of the last
$n$ matchings, the transitive closure of the "matched" relation defines two
sets $M^*$ and $N^*$ each having $m$ members. Each member in $M^*$ is a
subgraph containing the edges of an equivalence class, with exactly one edge
from each $E^*_r$ where $r$ is a descendent of $j$. Each member of $N^*$ is a
subgraph containing the edges of an equivalence class, with exactly one edge
from each $E^*_r$ where $r$ is {\em not} a descendent of $j$ (note that if
$j=1$ then $i=1$ and since $N(v,1,1)=0$ always, we may assume $j > 1$, and
thus $N^*$ is not empty). Note that the matchings $B_j(u)$ for all $u \in V$
match the members of $M^*$ with the members of $N^*$, and each such match
produces a member of $L^*$. Let us estimate $|N(v,i,j)|$ given that we know
exactly what $N^*$ contains; i.e. we shall estimate $\{|N(v,i,j)| \;\; | \;\;
N^*\}$. Consider a set $U=\{(x_1,u_1),(x_2,u_2),\ldots,(x_k,u_k)\}$ of $k$
edges, where for $t=1,\ldots,k$, $(x_t,u_t) \in D^-_{p(j)}(u_t)$, and
$(x_t,u_t)$ belongs to a member $T^t$ of $N^*$ containing an edge of
$D^+_i(v)$. The last requirement is valid since all the edges of $D^+_i(v)$
belong to members of $N^*$ because $i$ is not a descendent of $j$. Similarly,
the edges of $D^-_{p(j)}(u_t)$ belong to members of $N^*$ since $p(j)$ is not
a descendent of $j$. We call $U$ {\em bad}, if for all $t=1,\ldots,k$, $(x_t,u_t)$
is matched in $B_j(u_t)$ to an edge $(u_t,y_t) \in D^+_j(u_t)$ where $y_t$
already appears in $T^t$. (Note that the edges in $D^+_j(u_t)$ belong to
members of $M^*$). Since there are less than $h$ vertices in $T^t$, and since
$B_j(u_t)$ is selected at random, we have that
$$
Prob[(x_t,u_t) \;\; is \;\; matched \;\; in \;\; B_j(u_t) \;\; to \;\; a
\;\; bad \;\; edge] \leq \frac{h}{d^+_j(u_t)}.
$$
Similarly, the probability that $(x_t,u_t)$ is matched in $B_j(u_t)$ to a bad
edge, {\em given that} $(x_s,u_s)$ is matched in $B_j(u_s)$ to a bad edge,
for all $1 \leq s < t$, is at most $h/(d^+_j(u_t) - (t-1))$. Thus,
$$
{\rm Prob}[U \;\; is \;\; bad] < \Pi_{t=1}^{k} \frac{h}{d^+_j(u_t)-t+1}.
$$
Assuming $k \leq d^+_j(u_t)/2$, and using (\ref{e6}) we have
$$
{\rm Prob}[U \;\; is \;\; bad] < \left(\frac{12h^2}{C_1 \log n}\right)^k.
$$
Consequently,
$$
{\rm Prob}[|N(v,i,j)| \geq k \;\; | \;\; N^*] < {d^+_i(v) \choose
k}\left(\frac{12h^2}{C_1 \log n}\right)^k.
$$
Note that the estimation in the last inequality does not depend on $N^*$, and
thus,
$$
\Pr [|N(v,i,j)| \geq k] < {d^+_i(v) \choose
k}\left(\frac{12h^2}{C_1 \log n}\right)^k.
$$
Now put $k=\lfloor \log n \rfloor$ and note that indeed,
$k \leq d^+_j(u_t)/2$. We may therefore estimate the last inequality as follows:
$$
\Pr [|N(v,i,j)| \geq \log n] < 2^{d^+_i(v)}
\left(\frac{12h^2}{C_1 \log n}\right)^{\log n} < 2^{C_1\log n}\left(\frac{12h^2}{C_1 \log n}\right)^{\log n}=
$$
$$
\left(\frac{2^{C_1} \cdot 12h^2}{C_1 \log n}\right)^{\log n} < \frac{1}{10nh^2}.
$$
Thus, with probability at least $1 - nh^2/(10nh^2) \geq 0.9$, for all $1 \leq
i \leq j \leq h$, and for all $v \in V$, $|N(v,i,j)| \leq \log n$. \npf

For two vertices $u,v \in V$ (not necessarily distinct) and for two indices
$0 \leq j < i \leq h$ let $L([u,j],[v,i])$ denote the set of members of
$L^*$ which contain an edge of $D^-_i(v)$ and also contain an edge of
$D^-_j(u)$. Note that when $j=0$, $D^-_j(u)$ is undefined, so we define
$D^-_0(u)=D^+_1(u)$ and $d^-_0(u)=d^+_1(u)$ in this case only. For the sake
of symmetry, define $L([v,i],[u,j])=L([u,j],[v,i])$, and define
$L([u,i],[v,i]) = 0$, when $u \neq v$.

\begin{lemma}
\label{la4}
If all the perfect matchings $B_i(v)$ are selected randomly and
independently, then, with probability at least 3/4, for every $u,v \in V$ and
for $0 \leq j < i \leq h$,
\begin{equation}
\label{e7}
|L([u,j],[v,i])| \leq \frac{C_1}{10h^5}\log n.
\end{equation}
\end{lemma}
{\bf Proof:}\, Consider first the case where $j=p(i)$ or $i=1$ and $j=0$. In
this case, $L([u,j],[v,i])$ is simply the set of members of $L^*$ which
contain $(u,v)$ as their edge from $E^*_i$. Trivially, this set is empty
if $(u,v) \notin E^*_i$ and contains exactly one element if $(u,v) \in
E^*_i$. Thus, $|L([u,j],[v,i])| \leq 1$ in this case, so (\ref{e7}) clearly
holds.

\noindent
We may now assume $i > 1$ and $j \neq p(i)$. Let $k=p(i)$, so we must have $k
\neq j$. Suppose that we know, {\em for all $x \in V$}, that
$|L([u,j],[x,k])|=f_x$ (i.e. we know all these $n$ values). We wish to
estimate the value of $|L([u,j],[v,i])|$ given this knowledge. This is done
as follows. Let $L([u,j],[x,k],[v,i])$ be the subset of $L([u,j],[v,i])$
consisting of the members having an edge of $D^-_k(x)$. Note that
$|L([u,j],[x,k],[v,i])| \leq |L([x,k],[v,i])| \leq 1$ according to the
previous case, since $k=p(i)$. More precisely, if $(x,v) \notin D^-_i(v)$
then $|L([u,j],[x,k],[v,i])|=0$. If, however, $(x,v) \in D^-_i(v)$ then
$$
E[|L([u,j],[x,k],[v,i])| \;\;\;\; | \;\;\;\; |L([u,j],[x,k])| = f_x ] = f_x /
d^+_i(x),
$$
since the matching $B_i(x)$ is selected at random and $f_x/d^+_i(x)$ is the
probability that $(x,v)$ is matched to one of the $f_x$ members of $D^-_k(x)$
which are edges of members of $L([u,j],[x,k])$. Thus, if we put
$$
R_x = \{ |L([u,j],[x,k],[v,i])| \;\;\;\; | \;\;\;\; |L([u,j],[x,k])| = f_x \}
$$
then for $(x,v) \in D^-_i(v)$ we have that $R_x$ is an indicator random
variable with $E[R_x] = {\rm Prob}[R_x = 1] = f_x/d^+_i(x)$, while for $(x,v)
\notin D^-_i(v)$ we have $R_x = 0$. Note that if $(x,v) \in D^-_i(v)$ and $x
\neq y$ then $R_x$ is independent from $R_y$, since the value of $R_x$
depends only on the matching $B_i(x)$, which is independent from the matching
$B_i(y)$. Let
$$
R = \{ |L([u,j],[v,i])| \;\;\;\; | \;\;\;\; \forall x \in V,
|L([u,j],[x,k])|=f_x \}.
$$
According to the definition of $R$, we have
$$
R = \sum_{x \in V} R_x = \sum_{(x,v) \in D^-_i(v)} R_x.
$$
Thus, $R$ is the sum of independent indicator random variables. By linearity
of expectation,
$$
E[R] = \sum_{(x,v) \in D^-_i(v)} E[R_x] = \sum_{(x,v) \in D^-_i(v)}
f_x/d^+_i(x).
$$
On the other hand, we know that $\sum_{x \in V}f_x = d^-_j(u)$, since this
sum equals to the number of copies of $L^*$ having an edge of $D^-_j(u)$, and
this number is exactly $d^-_j(u)$. We also know from (\ref{e6}) that
$d^+_i(x) \geq 0.17\frac{C_1}{h}\log n$. Therefore,
\begin{equation}
\label{e8}
E[R] \leq \frac{d^-_j(u) h}{0.17C_1\log n}.
\end{equation}
Note that if $f_x = 0$ for some $x \in V$, then $R_x = 0$, and the term $R_x$
can be eliminated from the sum which yields $R$. Since $\sum_{x \in V}f_x =
d^-_j(u)$ this means that $R$ is the sum of at most $d^-_j(u)$ independent
indicator random variables. We can now apply the Chernoff bounds for $R$, and
obtain, for every $\alpha > 0$:
$$
{\rm Prob}[R - E[R] > \alpha] < \exp\left(-\frac{2\alpha^2}{d^-_j(u)}\right).
$$
In particular, for $\alpha = \sqrt{d^-_j(u)\log(2hn)}$,
$$
{\rm Prob}[R - E[R] > \sqrt{d^-_j(u)\log(2hn)}] <
\exp\left(-\frac{2d^-_j(u)\log(2hn)}{d^-_j(u)}\right) = \frac{1}{4h^2n^2},
$$
and it now follows from (\ref{e8}) that with probability at least
$1-1/(4h^2n^2)$,
$$
R \leq  \frac{d^-_j(u) h}{0.17C_1\log n} + \sqrt{d^-_j(u)\log(2hn)}.
$$
However,
$$
d^-_j(u)=d_j(u)-d^+_j(u)\leq d_j(u)-0.17\frac{C_1}{h}\log n \leq
\frac{d(u)}{h}+0.05\frac{d(u)}{h^2}-0.17\frac{C_1}{h}\log n< 1.5\frac{C_1}{h}\log n
$$
and therefore, with probability at least $1-1/(4h^2n^2)$,
\begin{equation}
\label{e9}
R \leq 9+\sqrt{1.5\frac{C_1}{h}\log n\log(2hn)} \ll \frac{C_1}{10h^5}\log n.
\end{equation}
Note that the estimation for $R$ in (\ref{e9}) does not depend on the
$f_x$'s. Thus, with probability at least $1-1/(4h^2n^2)$,
$$
|L([u,j],[v,i])| \leq \frac{C_1}{10h^5}\log n.
$$
Consequently, with probability at least $1-h^2n^2/(4h^2n^2) = 3/4$,
(\ref{e7}) holds for all $u,v \in V$ and for $0 \leq j < i \leq h$.
\npf

\noindent
{\bf Completing the proof:}\, According to Lemmas \ref{la3} and
\ref{la4} we know that with probability at least 0.65, we can obtain a
decomposition $L^*$ with the properties guaranteed by Lemmas \ref{la3} and
\ref{la4}. We therefore fix such a decomposition, and denote it by $L'$. We
let each member $T \in L'$ choose an integer $c(T)$, where $1 \leq c(T) \leq
h$. Each value has equal probability $1/h$. All the $m$ choices are
independent. Let $C(v,i)$ be the set of members of $T$ which selected $i$ as
their value and they contain an edge of $D^+_i(v)$. Put $|C(v,i)|=c(v,i)$. Clearly,
$0 \leq c(v,i) \leq d^+_i(v)$, and $E[c(v,i)] = d^+_i(v)/h$. Since the
choices are independent, we know that
$$
{\rm Prob}\left[c(v,i) < \frac{d^+_i(v)}{h+1}\right] =
{\rm Prob}\left[c(v,i) -E[c(v,i)] < -\frac{d^+_i(v)}{h(h+1)}\right] <
$$
$$
\exp\left(-\frac{2d^+_i(v)^2}{(h+1)^2h^2d^+_i(v)}\right) <
\exp\left(-\frac{2d^+_i(v)}{(h+1)^4}\right) \leq
\exp\left(-\frac{0.17\frac{C_1}{h}\log n}{(h+1)^4}\right) <
\frac{1}{2nh}.
$$
Thus, with positive probability (in fact, with probability at least 0.5), we
have that for all $v \in V$ and for all $i=1,\ldots,h$,
\begin{equation}
\label{e10}
c(v,i) \geq \frac{d^+_i(v)}{h+1}.
\end{equation}
We therefore fix the choices $c(T)$ for all $T \in L'$ such that (\ref{e10})
holds.

\noindent
We are now ready to mend $L'$ into a decomposition $L$ consisting only of
trees. Recall that each member of $L'$ is homomorphic to $H(q)$. We shall
perform a process which, in each step, reduces the overall number of bad
edges in $L'$ by at least one. Thus, at the end, there will be no bad edges,
and all the members are, therefore, trees. Our process uses two sets $L_1$
and $L_2$ where, initially, $L_1 = L'$ and $L_2 = \emptyset$. We shall
maintain the invariant that, in each step in the process, $L_1 \cup L_2$ is a
decomposition of $G$ into subgraphs homomorphic to $H(q)$. Note that this
holds initially. We shall also maintain the property that $L_1 \subset L'$.
Our process halts when no member of $L_1 \cup L_2$ contains a bad edge, and
by putting $L = L_1 \cup L_2$ we obtain a decomposition of $G$ into copies of
$H$, as required. As long as there is a $T^\alpha \in L_1 \cup L_2$ which
contains a bad edge, we show how to select a member $T^\beta \in L_1$, and
how to create two subgraphs $T^\gamma$ and $T^\delta$ which are also
homomorphic to $H(q)$ with $E(T^\alpha) \cup E(T^\beta) = E(T^\gamma) \cup
E(T^\delta)$, such that the number of bad edges in $E(T^\gamma) \cup
E(T^\delta)$ is {\em less} than the number of bad edges in $E(T^\alpha) \cup
E(T^\beta)$. Thus, by deleting $T^\alpha$ and $T^\beta$ from $L_1 \cup L_2$
and inserting $T^\gamma$ and $T^\delta$ both into $L_2$, we see that $L_1
\cup L_2$ is a {\em better} decomposition since it has less bad edges. It
remains to show that this procedure can, indeed, be done.

\noindent
Let $i$ be the maximum number such that there exists a member $T^{\alpha} \in
L_1 \cup L_2$ where $T^{\alpha}(i)$ is bad. Let $T^{\alpha}(i)=(v,w)$.
Consider the subgraph $T^{\epsilon}$ of $T^{\alpha}$ consisting of all the
edges $T^{\alpha}(j)$ where $j$ is a descendent of $i$. Our aim is to find a
member $T^\beta \in L_1$, which satisfies the following requirements:
$$
1) ~ c(T^\beta)=i.
\qquad \qquad
2) ~ T^\beta(i) \in D^+_i(v).
\qquad \qquad
3) ~ {\rm No ~vertex ~ of ~} T^\alpha ~,~ {\rm except} ~ v ~,~ {\rm appears ~ in} ~ T^\beta.
$$
We show that such a $T^\beta$ can always be found. The set $C(v,i)$ is
exactly the set of members of $L'$ which meet the first two requirements
(although some of them may not be members of $L_1$). Let $U$ be the set of
vertices of $T^\alpha$, except $v$. For $u \in U$, and for all $0 \leq
j \leq h$, all the members of $L([u,j],[v,p(i)])$ are not allowed to be
candidates for $T^{\beta}$. This is because each member of
$L([u,j],[v,p(i)])$ contains an edge of $D^-_{p(i)}(v)$, and thus an edge of
$D^+_i(v)$, but it also contains the vertex $u$, which we want to avoid in
$T^\beta$, according to the third property required. According to Lemma
\ref{la4},
$$
|L([u,j],[v,p(i)])| \leq \frac{C_1}{10h^5}\log n.
$$
Hence,
$$
|\cup_{u \in U} \cup_{j=0}^{h-1} L([u,j],[v,p(i)])| < h^2\frac{C_1}{10h^5}{\log n}=\frac{C_1}{10h^3}\log n.
$$
Let $C'(v,i)$ be the set of members of $C(v,i)$ which satisfy the third
requirement. By (\ref{e10}), (\ref{e6}) and the last inequality,
$$
|C'(v,i)| \geq c(v,i) - \frac{C_1}{10h^3}\log n \geq \frac{d^+_i(v)}{h+1} -
\frac{C_1}{h^3}\log n \geq
$$
$$
0.17\frac{C_1}{h(h+1)}\log n - \frac{C_1}{10h^3}\log n > (h+1) \log n.
$$
We need to show that at least one of the members of $C'(v,i)$ is also in
$L_1$. Each member $T \in C(v,i)$ that was removed from $L'$ in a prior stage
was removed either because it had a bad edge $T(j)$ where $j \geq i$ (this is
due to the maximality of $i$), or because it was chosen as a $T^\beta$
counterpart of some prior $T^\alpha$, having a bad edge $T^\alpha(i)=(v,z)$
for some $z$. There are at most $\sum_{j=i}^h|N(v,i,j)|$ members $T \in C(v,i)$ which have
a bad edge $T(j)$ where $j \geq i$, and there are at most $|N(v,i,i)|$
members $T \in C(v,i)$ having $T(i)$ as a bad edge. According to Lemma
\ref{la3}, $|N(v,i,i)|+\sum_{j=i}^h |N(v,i,j)| \leq (h+1)\log n$. Since
$|C'(v,i)| > (h+1)\log n$, we have shown
that the desired $T^\beta$ can be selected.

\noindent
Let $T^{\pi}$ be the subgraph of $T^\beta$ consisting of all the edges
$T^{\beta}(j)$ where $j$ is a {\em descendent} of $i$. $T^\gamma$ is defined
by taking $T^\alpha$ and replacing its subgraph $T^\epsilon$ with the
subgraph $T^\pi$. Likewise, $T^\delta$ is defined by taking $T^\beta$ and
replacing its subgraph $T^\pi$ with the subgraph $T^\epsilon$. Note that
$T^\gamma$ and $T^\delta$ are both still homomorphic to $H(q)$, and that
$E(T^\alpha) \cup E(T^\beta) = E(T^\gamma) \cup E(T^\delta)$, so by deleting
$T^\alpha$ and $T^\beta$ from $L_1 \cup L_2$, and by inserting $T^\gamma$ and
$T^\delta$ to $L_2$ we have that $L_1 \cup L_2$ is still a valid
decomposition into subgraphs homomorphic to $H(q)$. The crucial point
however, is that every edge of $E(T^\alpha) \cup E(T^\beta)$ that was good,
remains good due to requirement 3 from $T^\beta$, and that the edge
$T^\alpha(i)$ which was bad, now plays the role of $T^\delta(i)$, and it is
now a good edge due to requirement 3. Thus, the overall number of bad edges
in $L_1 \cup L_2$ is reduced by at least one. \npf

\end{document}